\documentclass[10pt]{article} 
\usepackage{cancel}
\usepackage[english]{babel}
\usepackage[utf8x]{inputenc}
\usepackage{amsmath}
\usepackage{amssymb}
\usepackage{amsthm}
\usepackage{graphicx}
\usepackage[title]{appendix}
\usepackage{color}
\usepackage{amsfonts}
\usepackage{dsfont}
\usepackage{epstopdf}
\usepackage{float}
\usepackage{hyperref}
\usepackage[dvipsnames]{xcolor}
\usepackage{caption}
\usepackage{subcaption}
\usepackage{tikz}
\usepackage[pagewise]{lineno}

\usepackage{amssymb}

\newtheorem{defi}{Definition}[section]
\newtheorem{teo}[defi]{Theorem}
\newtheorem{prop}[defi]{Proposition}
\newtheorem{lem}[defi]{Lemma}
\newtheorem{cor}[defi]{Corollary}

\newtheorem{remark}[defi]{Remark}

\newcommand{\Om}{\Omega}
\newcommand{\OOO}{\Omega}
\newcommand{\ooo}{\overline}
\newcommand{\ep}{\varepsilon}
\newcommand{\AAAA}{\mathcal{A}}

\newcommand{\om}{\omega}
\newcommand{\dis}{\displaystyle}
\newcommand{\Fin}{\hfill$\Box$}
\newcommand{\R}{\mathds{R}}
\newcommand{\Supp}{\text{Supp\,}}

\newcommand{\Id}{\text{\rm Id.}}
\newcommand{\ppp}{\partial}
\begin{document}

\title{ \Large 
\textbf{Uniqueness in determining multidimensional domains with unknown initial data}}

\author{
J. Apraiz\thanks{Corresponding author. Universidad del Pa\'is Vasco, Facultad de Ciencia y Tecnolog\'ia, Dpto.\ Matem\'aticas, Barrio Sarriena s/n 48940 Leioa (Bizkaia), Spain. 
E-mail: {\tt jone.apraiz@ehu.eus}.},
\ \ A. Doubova\thanks{Universidad de Sevilla, Dpto.\ EDAN e IMUS, Campus Reina Mercedes, 41012~Sevilla, Spain, E-mail: {\tt doubova@us.es}.},
\ \ E. Fern\'andez-Cara\thanks{Universidad de Sevilla, Dpto.\ EDAN e IMUS, Campus Reina Mercedes, 41012~Sevilla, Spain, E-mail: {\tt cara@us.es}.}, 
\ \ M. Yamamoto\thanks{The University of Tokyo, Japan, E-mail:  {\tt myama@next.odn.ne.jp}.}
}

\date{}

\maketitle

\begin{abstract}
   This paper addresses several geometric inverse problems for some linear parabolic systems where the initial data (and sometimes also the coefficients of the equations) are unknown.
   The goal is to identify a subdomain within a multidimensional set.
   The non-homogeneous part of the equation is expressed as a function satisfying some specific assumptions near a positive time.
   We establish uniqueness results by incorporating observations that can be on a part of the boundary or in an interior (small) domain.
   Through this process, we also derive information about the initial data.
   The main tools required for the proofs include semigroup theory, unique continuation and time analyticity results.

 \end{abstract}

\vspace*{0,4in}
 


\section{Introduction}\label{Sec.Introd}

Let~$\Om\subset \mathds{R}^d$ be a $d$-dimensional bounded connected open set 
with smooth boundary ($d \ge 2$) and 
let~$D_1, D_2 \subset\subset \Om$ be two {\it simply connected} smooth 
subdomains. 

   For~$k = 1$ and~$k = 2$, we will consider the systems
   \begin{equation}\label{eq.pb}
\begin{cases}
\partial_t u_k + \mathcal{A} u_k = f(x,t) &  \text{in} \  (\Om\!\setminus\!
\overline{D}_k)\times(0,T), 
\\[1mm]
u_k = 0 & \text{on} \  \ppp D_k \times (0,T), \\[1mm]
u_k = g(x,t) & \text{on} \  \partial\OOO \times(0,T) ,
\end{cases}
   \end{equation}
where~$f \not\equiv 0$ can be viewed as an external source and $g \not\equiv 0$ is a boundary input.
   Both functions are assumed to activate the system with suitable designs. More specific forms of $f(x,t)$ and $g(x,t)$ will be considered in Theorems~\ref{th.main} and~\ref{th.main3} below and also in the comments at the end of Section \ref{Sec-proof1}.

   Note that the initial values of~$u_1$ and~$u_2$ are not specified. 

   In~\eqref{eq.pb}, we have used the following notation:
   \[
\mathcal{A} v(x) : = - \displaystyle\sum_{i,j=1}^d \partial_i(a_{ij}(x) 
\partial_j v(x)) 
+ 
\displaystyle\sum_{j=1}^d b_j(x) \partial_j v(x) + c(x) v(x) .
   \]
   We assume that the coefficients~$a_{ij} = a_{ji}$ belong to~$C^1(\overline{\Om})$, while the~$b_j$ and~$c$ belong to~$L^\infty(\Om)$. 
   Also, we assume that there exists~$\alpha > 0$ such that
   \begin{equation}\label{ellipticity}
\sum_{i,j=1}^d a_{ij}(x) \xi_i \xi_j 
\geq \alpha |\xi|^2 \quad \forall \xi \in \mathds{R}^d, \ \text{ a.e.
\ in~$\Om$}
   \end{equation}
and, at least, $f \in L^2(\Om \times (0,T))$ and~$g \in L^2(0,T;H^{3/2}(\partial\Om))$.
   
   For~$k = 1$ and~$k = 2$, we introduce the linear operators~$A_k :  \mathcal{D}(A_k) \to L^2(\Om\!\setminus\!\overline{D}_k)$, with
   \[
\mathcal{D}(A_k): = \{ v\in H_0^1(\Om\!\setminus\!\overline{D}_k) : 
\ \mathcal{A}v \in L^2(\Om\!\setminus\!\overline{D}_k) \}
   \]
and
   \[
(A_k v)(x) := \mathcal{A} v(x) \ \text{ a.e.\ in~$\Om\!\setminus\!
\overline{D}_k$ } \ \forall v\in \mathcal{D}(A_k) .
   \]

   In the sequel, we will denote by~$\frac{\partial}{\partial \nu_{A}}$ 
the conormal derivative associated to the coefficients~$a_{ij}$ 
on~$\partial\Om$, that is,
   \[
\frac{\partial v}{\partial \nu_A} := \sum_{i,j=1}^d a_{ij}\partial_i v \, \nu_j
   \]
   ($\nu = \nu(x)$ denotes the outward unit normal vector at points~$x \in 
\partial\Om$).
   Moreover, in the results that follow we will assume without loss of 
generality that
   \begin{equation}\label{3p}
c(x) \geq c_0>0 \ \text{ a.e.\ in~$\Om$}
   \end{equation}
for a constant~$c_0$ sufficiently large to have that the~$A_k^{-1}$ exist and 
are compact. 

\

In this paper we will first deal with the following uniqueness question:
   
\

\noindent
\textbf{Question~Q1 (uniqueness, same coefficients): }
{\it Let~$\gamma\subset \partial\Om$ be a nonempty open subboundary.
   Let~$u_1$ and~$u_2$ be weak solutions to~\eqref{eq.pb} respectively 
corresponding to~$D_1$ and~$D_2$.
   This means that, for each~$k\in\{1,2\}$, $u_k$ belongs to~$L^2(0,T;H^1 (\Om\setminus\overline{D}_k))$, the lateral trace of $u_k$ on~$\partial\Om \times (0,T)$ is equal to~$g$ and~$u_k$} satisfies the parabolic equation in~\eqref{eq.pb} in the distributional sense.
   Then,
   \begin{equation*}
\text{Does the equality~$\dfrac{\partial u_1}{\partial \nu_A} = \dfrac{\partial
 u_2}{\partial \nu_A}$ on~$\gamma\times (0,T)$ imply~$D_1 = D_2$?}
   \end{equation*}

\

   We emphasize that our main purpose is to establish uniqueness for an inverse problem for a parabolic system where the initial values are unknown.
   The motivation is clear: it is often difficult in practice to a priori specify initial values.

   Note that any weak solution to~\eqref{eq.pb} satisfies
   \[
 \partial_t u_k \in L^2(0,T;H^{-1}(\Om \setminus \overline{D}_k))
   \]
and consequently can be viewed as a continuous~$L^2(\Om \setminus \overline{D}_k)$-valued function on~$[0,T]$.
   In particular, $u_k(\cdot\,,0) \in L^2(\Om \setminus \overline{D}_k)$.
   
   Moreover, $u_k \in L^2(\delta,T;H^{2}(\Om \setminus \overline{D}_k))$ for any small~$\delta > 0$.
   This allows to give a sense to the conormal derivative of~$u_k$ on~$\gamma\times (0,T)$, for instance as a function in~$L^2_{\rm loc}(0,T;L^2(\gamma))$.

\

   In order to formulate a second uniqueness question, we have to work  with two {\it a priori} different partial differential operators in space. Thus, let the coefficients $a_{ij}^k = a_{ji}^k \in C^1(\overline{\Om})$, 
$b_j^k \in L^\infty(\Om)$ and~$c^k \in L^\infty(\Om)$ be given for~$k = 1,2$, 
with the~$a_{ij}^k$ satisfying~\eqref{ellipticity} and the~$c^k$ satisfying
~\eqref{3p}.
   Let us set
   \begin{equation*}
\mathcal{A}^k v(x) : = - \displaystyle\sum_{i,j=1}^d 
\partial_i(a_{ij}^k(x) \partial_j v(x)) 
+ \displaystyle\sum_{j=1}^d b^k_j(x) \partial_j v(x) + c^k(x) v(x)
    \end{equation*}
and let us introduce the operators~$P_k : \mathcal{D}(P_k) \to L^2(\Om\!\setminus\!\overline{D}_k)$ as before:
   \[
\mathcal{D}(P_k): = \{ v\in H_0^1(\Om\!\setminus\!\overline{D}_k) : \ 
\mathcal{A}^k v \in L^2(\Om\!\setminus\!\overline{D}_k) \}
   \]
and
   \[
(P_k v)(x) := \mathcal{A}^k v(x) \ 
\text{ a.e.\ in~$\Om\!\setminus\!\overline{D}_k$ }
   \]
for all~$v\in \mathcal{D}(P_k)$.

   Again, it will be assumed that $c_0$ is large enough so that the~$P_k^{-1}$ exist and are compact.
   
   Our second question is the following:

\

\noindent
\textbf{Question~Q2 (uniqueness, different coefficients): }{\it
   Let the set~$\om \subset\subset \Om\setminus(\overline{D_1 \cup D_2})$ be nonempty and open.
   Let~$u_k$ be a weak solution to
   \begin{equation}\label{3a}
\begin{cases}
\partial_t u_k + \mathcal{A}^k u_k = f(x,t) &  \text{in} \  (\Om\!\setminus\!
\overline{D}_k)\times(0,T), 
\\[1mm]
u_k = 0 & \text{on} \  \ppp D_k \times (0,T),
\\[1mm]
 u_k = g(x,t) & \text{on} \  \partial\OOO \times(0,T)
\end{cases}
   \end{equation}
for~$k = 1, 2$.
   Then,
   \begin{equation*}
\text{Does the equality~$u_1 = u_2$ in~$\om\!\times\!(0,T)$ imply
~$D_1 = D_2$?}
   \end{equation*}
   }

\

   Inverse problems of these kinds have been considered in the literature since several decades;
   see for instance~\cite{Ale_1}--\cite{Apraiz} and~\cite{Paternain, Yamamoto-1, Yamamoto-2}.
   Many of them are still awaiting for theoretical and numerical research.

   Note that, in the previous works, uniqueness questions are usually raised for solutions that satisfy prescribed initial conditions, contrarily to what we do here.
   In fact, in~\textbf{Q2}, we can interpret that even the equations satisfied by~$u_1$ and~$u_2$ are unknown.

   This paper is organized as follows. First, in Section~\ref{Sec.Main}, 
we present our main results:
   positive answers for~\textbf{Q1} and~\textbf{Q2} under various assumptions.
   Then, Sections~\ref{Sec-proof1} and~\ref{Sec-proof3} contain the proofs and several related complements.   
   Finally, in~Section~\ref{Exten} we have included some related extensions, 
comments and open problems.

\

\section{Main results}\label{Sec.Main}

   In the main results of this paper, the following assumptions will be 
imposed to~$f$:
   \begin{equation}\label{7a}
\text{$f, \partial_t f, \dots, \partial_t^m f \in L^2(\Om \times (0,T))$ for 
some~$m \geq 0$}
   \end{equation}
and~$\partial_t^m f$ is piecewise constant in time and discontinuous near some~$t_1$, that is,

\begin{equation}\label{7b}
\left\{
\begin{array}{l}
\text{There exist $t_0, t_1, t_2$ with $0 < t_0 < t_1 < t_2 \leq T$ such that}
\\ \noalign{\smallskip}
\partial_t^m f(x,t) \!=\! 
\left\{ \begin{array}{l}
a_1f_0(x) + r_1(x,t) \ \text{ for } \ t_0 < t < t_1, \\
a_2f_0(x) + r_2(x,t) \ \text{ for } \ t_1 < t < t_2, \\
\end{array}\right.
\\ \noalign{\smallskip}
\text{where $r_1: (t_0,t_1] \to L^2(\OOO)$ is analytic,
$r_2 \in L^2(\OOO \times (t_1,t_2))$,}
\\ \noalign{\smallskip}
a_1, a_2 \in \R, \ f_0 \in L^2(\Om) \ \text{ and } \ a_1f_0(x) \not\equiv a_2f_0(x) .
\end{array}\right.
\end{equation}

   Here, the assumptions on~$r_1$ is understood in the sense that there exists~$\ep > 0$ such that~$r_1$ can be extended to an analytical function in~$(t_0,t_1 + \ep)$.
   Also, note that~$f_0$ cannot be identically zero and~$a_1$ cannot be equal to~$a_2$.
   
   For simplicity, we will assume in this section that~$g$ is independent of~$t$ and, of course, $g \in H^{3/2}(\partial\Om)$.

\

Our first result furnishes an answer to~\textbf{Q1}.
   It is the following:

\

\begin{teo}\label{th.main}
   Let~$u_1$ and~$u_2$ be solutions to~\eqref{eq.pb} respectively 
corresponding to the simply connected open sets~$D_1$ and~$D_2$.
   Suppose that~$f$ satisfies~\eqref{7a}, \eqref{7b} and moreover the functions~$f_0$, $r_1$ and~$r_2$ in~\eqref{7b} satisfy
\begin{equation}\label{7c}
\left\{
\begin{array}{l}
\text{$f_0(x) = 0$ in~$D_1 \cup D_2$, $r_1(x,t) = 0$ in $(D_1 \cup D_2) \times (t_0,t_1)$}
\\ \noalign{\smallskip}
\text{and $r_2(x,t) = 0$ in $(D_1 \cup D_2) \times (t_1,t_2)$.}
\end{array}
\right.
   \end{equation}
   Then the answer to~${\bf Q1}$ is yes.
   Moreover, $u_1(\cdot\,,0) = u_2(\cdot\,,0)$ in $\OOO 
\setminus (\ooo{D_1\cup D_2})$.
\end{teo}

   In~Theorem~\ref{th.main} we assume that $\partial_t^m f$ possesses a jump discontinuity in time.
   For example, when~$m = 0$, $r_1$ and~$r_2$ vanish identically and $a_1=1$ and~$a_2 = 0$, we can interpret that~$f$ is a {\it bang-bang or on-off} control input, which is strongly motivated from a pratical viewpoint.

   It will be shown at the end of Section 3 that uniqueness may fail if for instance~$g \equiv 0$ and~$f$ only satisfies~\eqref{7a} and
   \[
f(x,t) = 0 \ \text{ in } \ (D_1 \cup D_2) \times (T_1,T_2)
   \]
for some~$T_1$ and~$T_2$ with $0 \leq T_1 < T_2 \leq T$.
   Hence, an additional condition on~$f$ (like~\eqref{7b} or something similar) is needed in the present context.

Actually, from the viewpoint of applications, it is acceptable to assume that we can choose~$f$ and~$g$ to perform adequate tests. We can take for example~$f(x,t) = f_0(x) \mu(t)$, with a structure of separation of variables, for some $f_0$ and~$\mu$.
   Note that inputs of this kind are very common (and useful) in the analysis of control problems, see for instance~\cite{Coron}.
   As shown below, they also seem reasonable for the analysis of inverse problems.
   
   On the other hand, it will be seen that some property on~$f$ apart from~\eqref{7a} and~\eqref{7c} is needed.
   See to this respect the comments in the final part of~Section~\ref{Sec-proof1}.

\

\begin{remark}\label{rem_2.1}{\rm
   It may seem that the assumption~\eqref{7c} is not realistic, since it deals with the unknown sets~$D_1$ and~$D_2$.
   However, it may be viewed as a reasonable hypothesis on~$f$ if we admit to have some {\it a priori} information on the locations of~$D_1$ and~$D_2$.
   Thus, if~$\Om_0 \subset\subset \Om$ is a non-empty open set, we can replace~\eqref{7c} by
   \[
\text{$D_1, D_2 \subset \Om_0$ and $f(x,t) = 0$ in~$\Om_0 \times (t_0,t_2)$.}  
   \]
   This means that we know that~$D_1$ and~$D_2$ cannot be too close to the outer boundary $\ppp\OOO$.
   In other words, our formulation aims at the detection of cavities existing in a deep part of~$\OOO$.
   Note that~$f$ can be viewed as an external source intended to activate the system for the sake of identifying an unknown subdomain.
   Such an input should not disturb the original system and it is thus quite practical that the activating region does not touch $\overline{D_1 \cup D_2}$ (otherwise, our formulation would correspond to an {\it invasive} test, that should be avoided in practice).  
\Fin
}
\end{remark}

\

   A direct consequence of~Theorem~\ref{th.main} is the following:

\

\begin{cor}\label{th.main2}
   Let~$u_1$ and~$u_2$ be as in~Theorem~\ref{th.main}.
   Assume that~$f$ is of the form~$f(x,t) = f_0(x) \mu(t)$, where
   \[
\left\{
\begin{array}{l}\dis
\text{$f_0 \in L^2(\Om)$, $f_0(x) = 0$ in~$D_1 \cup D_2$ and~$f_0 \not\equiv 
0$,}
\\ \noalign{\smallskip} \dis
\text{$\mu$ is piecewise polynomial and~$\mu \not\in C^m([0,T])$ for some~$m 
\geq 0$.}
\end{array}
\right.
   \]
   Then, again, the answer to~${\bf Q1}$ is yes and~$u_1(\cdot\,,0) 
= u_2(\cdot\,,0)$.
\end{cor}

   Here, that~$\mu$ is ``piecewise polynomial'' means that there exists a 
partition~$0 < t_0 < t_1 < \cdots < t_N = T$ such that~$\mu|_{(t_{i-1},t_i)}$ 
is polynomial for each~$i\in\{1,\ldots,N\}$.

\

  For~\textbf{Q2}, the following holds:

\

\begin{teo}\label{th.main3}
   Let~$\om$, $D_1$ and~$D_2$ be as above and let~$u_k$ be a weak solution 
to~\eqref{3a} for~$k = 1, 2$.
   Assume that~$f$ satisfies~\eqref{7a}, \eqref{7b} and moreover the functions~$f_0$, $r_1$ and~$r_2$ in~\eqref{7b} satisfy
  \begin{equation}\label{7d}
\left\{
\begin{array}{l}
\text{$f_0(x) = 0$ in~$D_1 \cup D_2 \cup \om$, $r_1(x,t) = 0$ in $(D_1 \cup D_2 \cup \om) \times (t_0,t_1)$}
\\ \noalign{\smallskip}
\text{and $r_2(x,t) = 0$ in $(D_1 \cup D_2 \cup \om) \times (t_1,t_2)$.}
\end{array}
\right.
   \end{equation}
   Also, assume that
   \begin{equation}\label{4pp}
P_1 P_2 v = P_2 P_1 v \quad \forall v \in C_0^\infty(\Om \setminus 
(\overline{D_1 \cup D_2})) .
   \end{equation}
   Then the answer to~${\bf Q2}$ is yes.
\end{teo}

   Thus, in accordance with Theorem~\ref{th.main3}, the uniqueness of~$D$ holds also when the elliptic operators~$\AAAA^1$ and~$\AAAA^2$ are unknown.
   In other words, the values taken by the solution in~$\om \times (0,T)$ contain more information on the geometry of~$D$ than the coefficients.
 
   Note that~\eqref{4pp} is trivially satisfied in the particular case of 
constant coefficients~$a_{ij}^k$, $b_i^k$ and~$c^k$ for~$k=1,2$ and~$1 \leq i,j \leq d$.

\

\section{Proof of Theorem~\ref{th.main}}\label{Sec-proof1}

   We will assume that~$D_1\neq D_2$ and we will prove the result by 
contradiction.

\

\textsc{Step 1:}
   Let us set~$G := \Om\setminus(\overline{D_1\cup D_2})$ and~$u:=u_1-u_2$ 
in~$G\times (0,T)$.
   Then one has
   \[
\begin{cases}
\partial_t u + \mathcal{A} u  = 0 & 	\text{in }G\times(0,T), 
\\[1mm]
u = 0                                  & \text{on }\partial \Om \times(0,T),
\\[1mm]
\dfrac{\partial u}{\partial \nu_A} = 0  & \text{on }  \gamma \times(0,T).
\end{cases}
   \]
 
   The {\it unique continuation principle} 
(see for instance~\cite{SautScheurer}) yields~$u = 0$ in~$G\times(0,T)$ and 
therefore~$u_1 = u_2$ in~$G\times(0,T)$. 
 
   Since~$D_1 \not= D_2$, there must exist~$k$ and a nonempty open set~$E 
\subset (\Om \setminus\overline{D}_k)\cap D_{3-k}$ such that~$\partial E 
\subset \partial D_k \cup \partial D_{3-k}$.
   In particular, if~$D_1 \cap D_2 = \emptyset$, we can take~$E = D_{3-k}$.
   
\begin{figure}
 \centering
 \includegraphics[height = 8cm]{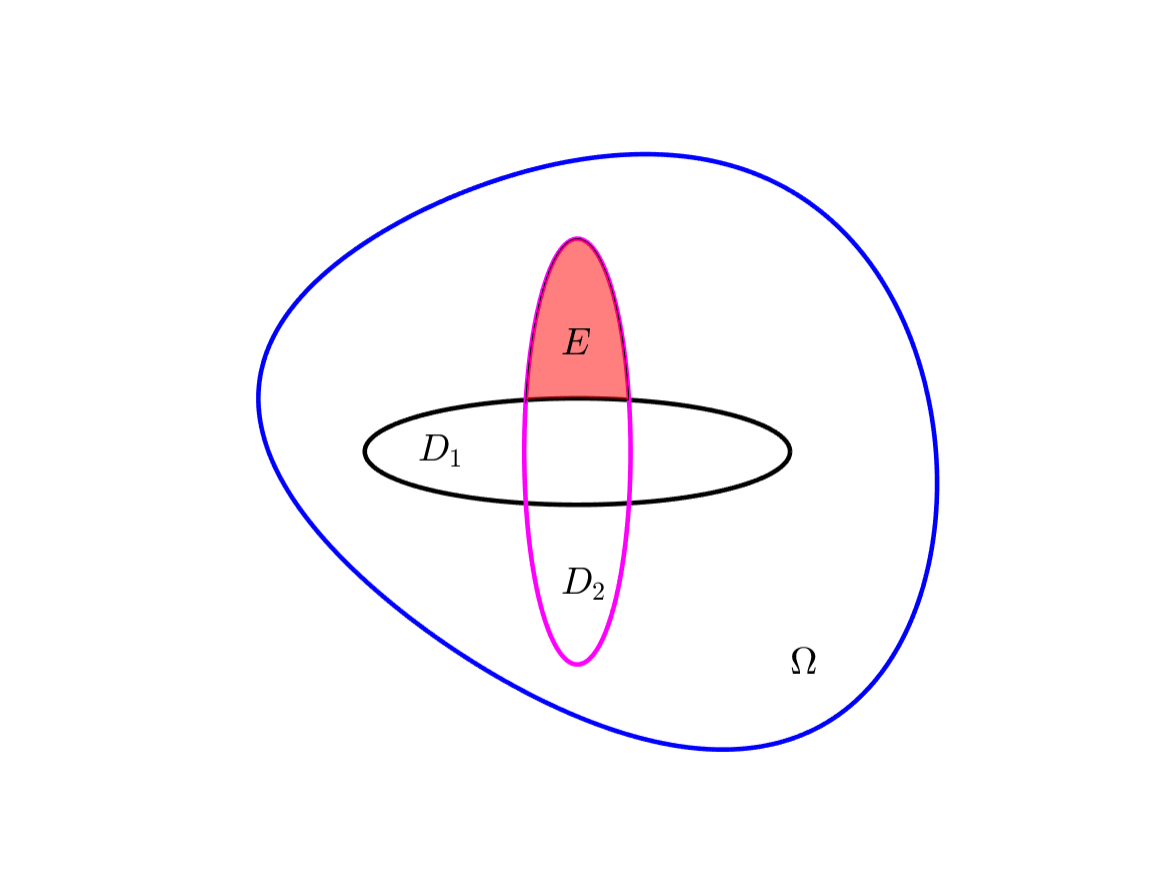}
 \caption{The set~$E$ with~$k=1$. }
 \label{fig.DomainE}
 \end{figure}

   Without loss of generality, we can asume that~$k = 1$, that is,
   \[
E \subset (\Om \setminus \overline{D}_1) \cap D_2 \ \text{ and } \ \partial E 
\subset \partial D_1 \cup \partial D_2;
   \]
   see~Figure~\ref{fig.DomainE}.


   Since~$u_1 = u_2$ in~$G$ and~$u_2 = 0$ on~$\partial D_2 \times (0,T)$, we  know that
   \[
\begin{cases}
\partial_t u_1 + A_1 u_1 = f(x,t) &  \text{in} \  
(\Om\!\setminus\!\overline{D}_1)\times(0,T), 
\\[1mm]
u_1 = 0 & \text{on} \  \partial(\Om\!\setminus\!\overline{D}_1) \times(0,T).
\end{cases}
   \]
   
   Moreover, from the classical parabolic smoothing property, we also have
   \[
\left\{
\begin{array}{l} \dis
\partial_t^n u_1 \in L^2(\delta,T; H_0^1(\Om\setminus\overline{D}_1)) \cap 
C^0([\delta,T];L^2(\Om\setminus\overline{D}_1)) \ \text{ and } 
\\ \noalign{\smallskip} \dis
\partial_t^{n+1} u_1 \in L^2((\Om\setminus\overline{D}_1) \times (\delta,T))
\end{array}
\right.
   \]
for any small~$\delta > 0$ and any~$n = 0,1,\dots,m$;
   see for instance~\cite{Amann, Friedman}.
   
   Let us set~$z = \partial_t^m u_1$.
   Then one has
   \[
\begin{cases}
\partial_t z + A_1 z  = \partial_t^m f(x,t) \ \text{ in } \ (\Om\setminus
\overline{D}_1)\times(0,T), 
\\[1mm]
z \in C^0((0,T];L^2(\Om\setminus\overline{D}_1)).
\end{cases}
   \]
   
   We recall that the operator $A_1$ is defined in a subspace of~$L^2(\OOO\setminus\ooo{D_1})$.
   In partcular, we see that
   \begin{equation}\label{1.10}
z_0 := z(\cdot\,,t_0) \in L^2(\Om\setminus\overline{D}_1)
   \end{equation}
and also
   \[
\begin{cases}
\partial_t z + A_1 z  = 0 & 	\text{in }E\times(t_0,t_2), 
\\[1mm]
z = 0 & \text{on }\partial E\times(t_0,t_2) 
\end{cases}
   \]
because $E \subset D_2$ and $f=0$ in $D_2 \times (0,T)$
by \eqref{7c}.

   Let~$E_0 \subset \subset E$ be a new nonempty open set, let us denote 
by~$\chi_0$ the restriction operator associated with~$E_0$, that is,
   \[
\chi_0 v : = v|_{E_0} \quad \forall v\in L^2(\Om\!\setminus\!\overline{D}_1)
   \]
and let us set
   \[
h(t) := \chi_0 z(\cdot\,,t) \quad \forall t \in (t_0,t_2).
   \]
   Then, from semigroup theory (see for instance~Proposition~2.1.9 
in~\cite{Lunardi}), we know that the~$L^2(E_0)$-valued function~
$t \mapsto h(t)$ is analytical in~$(t_0,t_2)$.
   
\

\textsc{Step 2:}
   In this step, we will prove that~$f_0 \equiv 0$, which is in contradiction with~\eqref{7b}
   (recall that a consequence of the last line of~\eqref{7b} is that~$f_0$ cannot vanish identically).
   
   After this, the proof of Theorem~\ref{th.main} will be complete.
   
   First, we recall two lemmas:

\

\begin{lem}\label{lem_3.1}
   Let us assume that~$r: (t_0,t_1) \to L^2(\Om \setminus \overline{D_1})$ possesses an analytic extension in~$(t_0,t_1+\ep)$ with~$\ep > 0$.
   Then
   \[
t \to \int^t_{t_0} e^{-(t-s)A_1} r(\cdot\,,s) \, ds
   \]
can also be extended to an analytic function in~$(t_0,\, t_1+\ep)$.
\end{lem}

\

   This result is an easy consequence of the fact that the semigroup~$\{ e^{-\tau A_1} \}$ is well-defined and holomorphic in a sector of the form
   \[
S_{\theta} := \{ \tau \in \mathbb{C} : \arg(\tau-t_0) < \theta, \ \vert \tau-t_0\vert < t_1-t_0+\ep \}
   \]
for some~$\theta \in (0,\pi/2)$.

\

\begin{lem}\label{lem_3.2}
   Let us assume that
   
\begin{enumerate}

\item $S = S(t)$ is real-valued and analytic in $t\in (t_0, \, t_1+\ep)$.

\item Its restriction $s := S\vert_{(t_0,t_1)}$ can be analytically extended to~$\tilde{s}$ in~$(t_0, t_1+\ep)$.

\end{enumerate}
   Then $\tilde{s}(t) = S(t)$ in~$(t_1, t_1+\ep)$.
\end{lem}

   Of course, Lemma~\ref{lem_3.2} is a direct consequence of the uniqueness of the holomorphic 
extension of $s$ from $(t_0,t_1)$ to~$(t_0, t_1+\ep)$. 

\

   Using again semigroup theory, \eqref{1.10}, \eqref{7b} and the definition of~$h$, we obtain the following identity in~$L^2(\Om \setminus \overline{D}_1)$:
   \[
 h(t) = \chi_0 e^{-(t-t_0) A_1} z_0 
 + \chi_0 \int_{t_0}^t e^{-(t-s)A_1}\partial_t^m f(\cdot\,,s) \, ds \quad 
\forall t \in (t_0,t_2) .
   \]

   Let us set 
   \[
r(x,t)=
\left\{
\begin{array}{ll}
r_1(x,t) & \text{ if } \ t_0<t<t_1,
\\
r_2(x,t) & \text{ if } \ t_1<t<t_2,
\end{array}
\right.
   \]
let $\zeta : (t_0,t_2) \to L^2(E_0)$ be given by
   \[
\zeta(t) = \chi_0 \int_{t_0}^t e^{-(t-s)A_1} r(\cdot\,,s) \, ds \quad \forall t \in (t_0,t_2)
   \]
and let us set
   \[
S(t) = h(t) - \chi_0 e^{-(t-t_0) A_1} z_0 - \zeta(t)
\quad \forall t \in (t_0,t_2) .
   \]
   
   Then, $\zeta$ is the restriction to~$E_0 \times (t_0,t_2)$ of the solution to
      \[
\begin{cases}
\partial_t Z + A_1 Z  = r(x,t) & 	\text{in }\Om\setminus\overline{D}_1\times(t_0,t_2), 
\\[1mm]
Z = 0 & \text{on }\partial(\Om\setminus\overline{D}_1)\times(t_0,t_2),
\\[1mm]
Z(x,t_0) = 0 & \text{in }\Om\setminus\overline{D}_1 .
\end{cases}
   \]
   
   Consequently, since $r(x,t) = 0$ in $(D_1 \cup D_2) \times (t_0,t_2)$, we have that~$\zeta$ and~$S$ are analytical in~$(t_0,t_2)$.

   Let~$s_1$ and~$s_2$ be the restrictions of~$S$ respectively to~$(t_0,t_1)$ and~$(t_1,t_2)$.
 Then, in view of~\eqref{7b}, one has
   \begin{equation}\label{s1}
\begin{array}{l}\dis
s_1(t) = \chi_0 \int_{t_0}^t e^{-(t-s)A_1} (\partial_t^m f(\cdot\,,s) - r_1(\cdot\,,s) )\, ds
\\ \dis
\phantom{s_1(t)} = a_1 \chi_0 A_1^{-1} \big(\text{Id.} - e^{-(t-t_0)A_1}\big) f_0
\end{array}
   \end{equation}
for all $t \in (t_0,t_1)$ and
   \begin{equation}\label{s2}
\begin{array}{l}\dis
s_2(t) = \chi_0 \int_{t_0}^t e^{-(t-s)A_1} (\partial_t^m f(\cdot\,,s) - r(\cdot\,,s) )\, ds
\\ \dis
\phantom{s_2(t)} = a_1 \chi_0 A_1^{-1}e^{-(t-t_1)A_1} \big(\text{Id.} - e^{-(t_1-t_0)A_1}\big) f_0 
\\ \noalign{\medskip} \dis
\phantom{s_2(t)} + a_2 \chi_0 A_1^{-1} \big(\text{Id.} - e^{-(t-t_1)A_1}\big) f_0
\end{array}
   \end{equation}
for all $t \in (t_1,t_2)$.

   From~Lemma~\ref{lem_3.1} we know that the function
   \[
t \to \chi_0\int^t_{t_0} e^{-(t-s)A_1}r_1(\cdot\,,s) \, ds
   \]
can be extended holomorphically to~$(t_0, t_1+\ep)$ for some~$\ep > 0$.
   Therefore, taking into account the definitions of~$S$ and~$s_1$ and~Lemma~\ref{lem_3.2}, we see that~$s_1$ can also be extended analytically to~$(t_0, t_1+\ep)$, with an extension that coincides with~$s_2$ in~$(t_1, t_1+\ep)$.
   
   Thus, thanks to~\eqref{s1}, \eqref{s2} and the fact that $a_1 \not= a_2$, one has
   \[
\chi_0 A_1^{-1} \big(\text{Id.} - e^{-(t-t_1)A_1}\big) f_0 = 0
\quad \forall t \in (t_1, t_1 + \ep)
   \]
and, after time differentiation, we see that
   \[
\chi_0 e^{-(t-t_1)A_1} f_0 = 0
\quad \forall t \in (t_1, t_1 + \ep) .
   \]
   
   This suffices to conclude that~$f_0$ vanishes in~$\OOO\setminus \ooo{D_1}$, which on the other hand is impossible.
   Indeed, if we introduce~$w$ with $w(\cdot\,,t):= e^{-tA_1}f_0$ for all $t \in (0, \ep)$, we obviously have
   \[
\left\{
\begin{array}{l}
\ppp_t w + \AAAA_1w = 0 \quad \mbox{in $(\OOO\setminus \ooo{D_1})
\times (0, \ep)$},
\\
w = 0 \quad \mbox{in $E_0 \times (0, \ep)$}
\end{array}
\right.
   \]
and, from unique continuation, $w(x,t) =  0$ in~$(\OOO\setminus \ooo{D_1})\times (0, \ep)$.
   But then, letting $t \to 0^+$, we finally get $f_0(x) = 0$ in~$\OOO\setminus \ooo{D_1}$.

    This ends the proof of~Theorem~\ref{th.main}.

\

   Let us complete this section with some complements to Theorem~\ref{th.main}:
   
\

   $\bullet$ In~\eqref{eq.pb} and~\textbf{Q1}, it also makes sense to have~$f \equiv 0$ and~$g \not\equiv 0$. 
   Let us see that in this case, under some appropriate conditions on~$g$, uniqueness holds again.
   
   Thus, let the functions~$u_k$ satisfy
   \[
\begin{cases}
\partial_t u_k + \mathcal{A} u_k = 0 &  \text{in} \  (\Om\!\setminus\!
\overline{D}_k)\times(0,T), 
\\[1mm]
u_k = g(x,t) & \text{on} \  \ppp (\OOO\setminus \ooo{D_k}) \times (0,T)
\end{cases}
   \]
and let us assume (for example) that~$g(x,t) = g_0(x)\mu(t)$ for all~$(x,t) \in \partial\Om \times (0,T)$, where~$g_0 \in H^{3/2}(\partial\Om)$ and
   \[
\mu(t)=
\left\{
\begin{array}{ll}
a_1 t, \quad &\text{if } 0 < t < t_1, \\
a_2 (t-t_1) + a_1 t_1, \quad &\text{if } t_1 < t < T
\end{array}
\right.
   \]
for some~$a_1, a_2 \in \R$ with~$a_1 \ne a_2$ and some~$t_1$ with~$0 < t_1 < T$.

   Let~$D_1, D_2 \subset \OOO$ be simply connected open sets such that~$\overline{D_1 \cup D_2} \subset \OOO$.
   There exist a neighborhood~$\Om_0$ of~$\ppp\OOO$ and a {\it lift}~$\tilde{g} \in H^2(\OOO)$ of~$g$ such that
   \[
\Supp\tilde{g} \subset \Om_0 \subset \ooo{\OOO} \setminus (D_1 \cup D_2) \ \text{ and } \ \tilde{g}\not\equiv 0 \ \text{ in } \ \Om_0.
   \]
   In particular, we have
   \begin{equation}\label{12}
\tilde{g} = 0 \quad \mbox{in $D_1 \cup D_2$}.
   \end{equation}
   
   Let us set~$v_k(x,t) := u_k(x,t) - \mu(t)\tilde{g}(x)$ in~$(\OOO\setminus \ooo{D_k}) \times (0,T)$.
   Then
   \begin{equation}\label{13}
\begin{cases}
\partial_t v_k + \mathcal{A} v_k = f(x,t) &  \text{in} \  (\Om\!\setminus\!
\overline{D}_k)\times(0,T), 
\\[1mm]
v_k = 0 & \text{on} \  \ppp (\OOO\setminus \ooo{D_k}) \times (0,T) .
\end{cases}
   \end{equation}
   
   In~\eqref{13}, we have introduced~$f(x,t) := -\mu'(t)\tilde{g}(x) - \mu(t)\AAAA\tilde{g}(x)$, that is,
   \[
f(x,t) = 
\left\{ \begin{array}{ll}
-a_1\tilde{g}(x) - r_1(t)\AAAA \tilde{g}(x) \quad &\text{if } t_0<t<t_1, \\
-a_2\tilde{g}(x) - r_2(t)\AAAA \tilde{g}(x) \quad &\text{if } t_1<t<t_2,
\end{array}\right.
   \]
where~$r_1(t) := a_1 t$ for~$t_0<t<t_1$ and~$r_2(t) := a_2 (t-t_1)+a_1 t_1$ for~$t_1<t<t_2$.

   Using~\eqref{12}, we see that~$f=0$ in~$(D_1\cup D_2) \times (0,T)$.
   Moreover, in view of the definition of~$\mu$ and~\eqref{12}, the conditions~\eqref{7a}--\eqref{7c} are satisfied for~$m = 1$.
   
   Therefore, Theorem~\ref{th.main} can be applied and one has~$D_1 = D_2$ and~$u_1(\cdot\,,0) = u_2(\cdot\,,0)$ in~$\OOO \setminus (\ooo{D_1 \cup D_2})$.

\

$\bullet$ Suppose again that~$f \equiv 0$ and~$g \not\equiv 0$ and assume now that we have the following additional information:
   \begin{equation}\label{extrainfo}
u_k(x,0) = 0 \ \text{ in $\Om\setminus\overline{D_k}$ for~$k=1,2$.}
   \end{equation}
   
   Then, from unique continuation, we first deduce (as in the proof of Theorem~\ref{th.main}) that~$u_1 = u_2$ in~$(\OOO \setminus \ooo{D_1 \cup D_2}) \times (0,T)$.
   
   Assume that~$D_1 \not= D_2$ and, for instance, $E := D_2\setminus \ooo{D}_1 \not= \emptyset$;
   then~$u_1$ satisfies the PDE
   \[
\partial_t u_1 + A_1 u_1 = 0 \ \text{ in } \ E \times (0,T),
   \]
    the boundary condition
   \[
u_1 = 0 \ \text{ on } \ \partial E \times (0,T)
   \]
and vanishes at~$t = 0$.
   Consequently, $u_1(x,t) = 0$ for~$x \in E$ and~$0 < t < T$ and, again from unique continuation, $u_1(x,t) = 0$ in $(\OOO\setminus \ooo{D}_1) \times (0,T)$, which contradicts that~$g = u_1\vert_{\ppp\OOO\times (0,T)} \not\equiv 0$.
   
  In other words, if~$f\equiv 0$, $g \not\equiv 0$ and for~$k = 1$ and~$k = 2$ the solution~$u_k$ to the problem~\eqref{eq.pb} vanishes at~$t = 0$, uniqueness holds.

   However, if the initial condition satisfied by the~$u_k$ is known but different from~$0$, the uniqueness assertion~$D_1 = D_2$ may fail, as the following example shows:
   take~$d=2$ and consider the square domains 
   \[
\begin{array}{c}
\OOO:= (1/4 , 2) \times (1/4 , 2), \\
D_1:= (1/2 , 1) \times (1/2 , 1)
\ \text{ and } \ D_2:= (1/2 , 3/2) \times (1/2 , 3/2)
\end{array}
   \]
and the function
   \[
u(x_1,x_2,t) := e^{-8\pi^2 t} \, \sin (2\pi x_1)\, \sin (2\pi x_2).
   \]
   
   It is then easy to check that~$\ppp_t u - \Delta u =  0$ in~$(\OOO\setminus \ooo{D}_k) \times (0,T)$ and~$u=0$ on~$\ppp D_k \times (0,T)$ for both $k=1$ and $k=2$.
   
   Hence, uniqueness does not hold in general.

\

$\bullet$ With~$f \not\equiv 0$ satisfying only~\eqref{7a} and~\eqref{7c} and~$g \equiv 0$, uniqueness can fail even under the additional assumption~\eqref{extrainfo}.

   This can be seen as follows.
   Assume that~$D_1 \not= D_2$ and let~$f_0 \in C_0^{\infty}(\OOO)$ be given with $f_0\not\equiv 0$ and~$\Supp f_0 \subset \OOO \setminus \ooo{(D_1\cup D_2)}$.
   Let us set $u(x,t) := tf_0(x)$ and $f(x,t):= f_0(x) + t\AAAA f_0(x)$.
   Then
   \[
\left\{ \begin{array}{ll}
\ppp_tu + \AAAA u = f(x,t) \quad& \mbox{in $(\OOO\setminus 
\ooo{D}_k) \times (0,T)$}, \\
u=0 \quad& \mbox{on $\ppp(\OOO\setminus \ooo{D}_k) \times (0,T)$},\\
u(\cdot,0)  = 0 \quad& \mbox{in $\OOO\setminus \ooo{D}_k$}
\end{array}\right.
   \]
for $k=1,2$.
   Since~$f_0$ can be chosen such that~$f_0 + t\AAAA f_0$ does not vanish in~
$(\OOO\setminus \ooo{(D_1 \cup D_2)}) \times (0,T)$, we see that uniqueness fails.

\

$\bullet$ On the other hand, if we assume that~\eqref{extrainfo} holds, uniqueness can be established for right-hand sides~$f(x,t) = f_0(x) \mu(t)$ with a smooth~$\mu$.

\

   This is shown in the following result:

\

\begin{prop}\label{prop2.7}
   Let us assume that $f(x,t) = f_0(x) \mu(t)$ a.e.\ with
   \[
\left\{
\begin{array}{l}
f_0 \in L^2(\OOO), \ \ \Supp f_0 \subset \OOO\setminus \ooo{(D_1 
\cup D_2)} \ \text{ and } \ f_0 \not\equiv 0,
\\
 \mu\in C^1([0,T]) \ \text{ and } \ \mu \not\equiv 0
 \end{array}
 \right.
    \]
  and~$g \equiv 0$.
    Also, let us assume that~\eqref{extrainfo} holds.
   Then, the answer to {\bf Q1} is yes.
\end{prop}

{\sc Proof:}
   Arguing as in the proof of~Theorem~\ref{th.main}, we see that~$u_1 = u_2$ in~$(\Om \setminus \overline{D_1 \cup D_2}) \times (0,T)$.
   
   Let us assume that~$D_1 \not= D_2$ and there exists a non-empty open set~$E \subset (\OOO\setminus \ooo{D}_1) \cap D_2$.
   Then~$f_0(x) = 0$ in $E$ and, since~$u_1(x,0)=0$ in $E$, one has
   \begin{equation}\label{n23}
u_1(x,t) = 0 \quad \mbox{in $E \times (0,T)$}.
   \end{equation}
   
   Let us consider the solution~$v$ to 
   \[
\left\{ \begin{array}{ll}
\ppp_tv + A_1v = 0 \quad& \mbox{in $(\OOO\setminus \ooo{D}_1)\times (0,T)$}, 
                                   \\
v = 0 \quad& \mbox{on $\ppp(\OOO\setminus \ooo{D}_1) \times (0,T)$},\\
v(\cdot\,,0) = f_0 \quad& \mbox{in $\OOO\setminus \ooo{D}_1$}.
\end{array}\right.
    \]
    Then, we can directly verify that
    \[
u_1(x,t) = \int^t_0 \mu(s)v(x,t-s) \, ds \ \text{ in } \ (\OOO\setminus \ooo{D}_1) \times (0,T).
   \]

   In view of~\eqref{n23}, for any~$\psi \in C^{\infty}_0(E)$ we have
   \[
\int^t_0 \mu(s)(v(\cdot\,,t-s),\psi)_{L^2(E)} = 0 \quad \forall t \in (0,T)
   \]
and, since the function~$s \mapsto (v(\cdot\,,s),\psi)_{L^2(E)}$ belongs to~$L^1(0,T)$ and $\mu \not\equiv 0$ in $(0,T)$, {\it Titchmarsh's Theorem} on convolutions (see~\cite{Ti}) implies
   \[
(v(\cdot\,,s)\,, \psi)_{L^2(E)} = 0\quad \forall s \in (0,T).
   \]  
But the function~$\psi$ is arbitrary in~$C^{\infty}_0(E)$. Consequently, $v(x,t) = 0$ in~$E \times (0,T)$ and (once more) from unique continuation, $v(x,t) = 0$ in~$(\OOO\setminus \ooo{D}_1)\times (0,T)$. Obviously, this implies~$f_0(x) =0$ in $\OOO\setminus \ooo{D}_1$, which is impossible. Thus, the proof of Proposition 2.7 is complete.
\Fin

\

   For the considered inverse problem, it also makes sense to try to retrieve~$D$ from the conormal derivative on a part of the lateral boundary.
   The related question is as follows:

\
   
\noindent
\textbf{Question~Q3 (reconstruction, same coefficients): }{\it
   Let~$\gamma\subset \partial\Om$ be a nonempty open subboundary and let~$u$ 
be a weak solution to 
   \begin{equation}\label{IP-2}
\begin{cases}
\partial_t u + \mathcal{A} u = f(x,t) &  \text{in} \  (\Om\!\setminus\!
\overline{D})\times(0,T), \\[1mm]
u = 0 & \text{on } \  \partial(\Om\!\setminus\!\overline{D}) \times(0,T),
\end{cases}
   \end{equation}
for some nonempty simply connected open set~$D \subset\subset \Om$.
   Assume that
   \begin{equation}\label{8p}
\frac{\partial u}{\partial \nu_A} = \beta \ \text{ on } \ \gamma\!\times\!
(0,T).
   \end{equation}
   Then,
   \[
\text{Can we find~$D$ (and~$u|_{t=0}$) from~$f$ and~$\beta$?}
   \]
   }
   
\

   From Theorem~\ref{th.main}, we see that an appropriate reformulation of 
the reconstruction problem may be the following:
   \[
\begin{cases}\dis
\text{Minimize } \frac{1}{2}
\left\Vert \frac{\partial u}{\partial \nu_A} - \beta \right\Vert_X^2
\\[3mm] \dis
\text{Subject to } D \in \mathcal{B}, 
\ u_0 \in L^2(\Om\setminus\overline{D}), \ u \ \text{ solves~\eqref{new},}
\end{cases}
   \]
where~$\beta$ is given, the admissible class of subdomains~$\mathcal{B}$ and 
the Hilbert space~$X$ are appropriately chosen and~\eqref{new} stands for the 
system
   \begin{equation}\label{new}
\begin{cases}
\partial_t u + \mathcal{A} u = f(x,t) &  \text{in} \  
\Om\!\setminus\!\overline{D}\times(0,T),
\\[1mm]
u = 0 & \text{on } \  \partial(\Om\!\setminus\!\overline{D}) \times(0,T),
\\[1mm]
u|_{t=0} = u_0 &  \text{in} \ \Om\!\setminus\!\overline{D} .
\end{cases}
   \end{equation}
   
   The numerical resolution of this extremal problem can be achieved by the 
method of fundamental solutions.
   This will be considered in a forthcoming paper.

Recall that this strategy has been applied in~\cite{Jairo} for the 
reconstruction of~$D$ in the context of a stationary inverse problem similar 
to~\eqref{IP-2}--\eqref{8p}.
   
\

\section{Proof of Theorem~\ref{th.main3}}\label{Sec-proof3}

   First, note that for any small~$\delta > 0$ one has
   \[
u_k, \partial_t u_k, \partial_i u_k, \partial_i\partial_j u_k \in L^2((\Om 
\setminus \overline{D}_k) \times (\delta,T))
   \]
for all~$k = 1,2$ and~$1 \leq i,j \leq d$, in view of the assumptions on~$f$.

   Let~$\om \subset \subset \Om\setminus(\overline{D_1 \cup D_2})$ be given.
   Without loss of generality, we can assume that there exists another 
nonempty open set~$G \subset \subset \Om\setminus(\overline{D_1 \cup D_2})$ 
such that~$\overline{\om} \cap \overline{G} = \emptyset$ and~$\Supp f 
\subset \overline{G} \times [0,T]$.
   
   Now, let~$\Om_0 \subset \Om$ be a new open set satisfying
   \[
\Om_0 \supset \overline{D_1 \cup D_2} , \quad \Om_0 \supset \overline{\om}, 
\quad \overline{\Om_0} \cap \overline{G} = \emptyset
   \]
and set~$O := \Om_0 \setminus (\overline{D_1 \cup D_2})$.
   Then, for~$k=1,2$,
   \begin{equation}\label{3.3}
\partial_t u_k + P_k u_k = 0 \ \text{ in } \ O \times (0,T),
   \end{equation}
where the~$P_k$ are the operators defined in Section~\ref{Sec.Introd}.

\

   We will use the following unique continuation result:

\begin{lem}\label{lem.5.1}
   Let~$T_1$ and~$T_2$ be given with~$T_1 < T_2$ and let us assume that 
the functions~$v$, $\partial_t v$, $\partial_i v$ and~$\partial_i\partial_j v$ 
belong to~$L^2(O \times (T_1,T_2))$.
   If
   \[
(\partial_t + P_1)(\partial_t + P_2)v = 0 \ \text{ in } \ O \times (T_1,T_2)
   \]
and~$v(x,t) = 0$ a.e.\ in~$\om \times (T_1,T_2)$, then~$v \equiv 0$.
\end{lem}

\noindent
{\sc Proof:}
   Set~$w := (\partial_t + P_2)v$.
   Then~$w \in L^2(O \times (T_1,T_2))$ and
   \[
\begin{cases}
\partial_t w + P_1 w = 0 &  \text{in} \  O\times(T_1,T_2), 
\\[1mm]
w = 0 & \text{in} \  \om\times(T_1,T_2) .
\end{cases}
   \]

   From the unique continuation property satisfied by~$\partial_t + P_1$, we must have~$w = 0$ in the whole horizontal component of~$\om \times (T_1,T_2)$, that is, in~$O \times (T_1,T_2)$.
   Hence,
   \[
\begin{cases}
\partial_t v + P_2 v = 0 &  \text{in} \  O\times(T_1,T_2), 
\\[1mm]
v = 0 & \text{in} \  \om\times(T_1,T_2)
\end{cases}
   \]
and, again from unique continuation, we obtain~$v = 0$ in~$O \times (T_1,T_2)$.
\Fin

\

   Let us end the proof of Theorem~\ref{th.main3}.
   
   From~\eqref{3.3} and the fact that~$P_1$ and~$P_2$ commute, we see that
   \[
(\partial_t + P_1)(\partial_t + P_2) = (\partial_t + P_2)(\partial_t + P_1)
   \]
and consequently
   \[
(\partial_t \!+\! P_1)(\partial_t \!+\! P_2)(u_1 \!-\! u_2)
= (\partial_t \!+\! P_2)(\partial_t \!+\! P_1)u_1 \!-\! (\partial_t \!+\! P_1)
(\partial_t \!+\! P_2)u_2 
= 0
   \]
in~$O \times (\delta,T)$ for any small~$\delta > 0$.
   Using Lemma~\ref{lem.5.1}, we find that~$u_1 = u_2$ in this set and, since 
$\delta$ is arbitrarily small, $u_1 = u_2$ in~$O \times (0,T)$.
   
   Let us assume that~$D_1 \not= D_2$.
   Then, without loss of generality, we can find again nonempty open sets
~$E_0$ and~$E$ satisfying
   \[
E_0 \subset \subset E \subset (O \setminus\overline{D}_1)\cap D_2 \ 
\text{ and } \ \partial E \subset \partial D_1 \cup \partial D_2 .
   \]
   
Arguing as we did in the case of Theorem~\ref{th.main}, we deduce that
   \[
\chi_0 \big(\text{Id.} - e^{-(t-t_1)P_1}\big) f_0 = 0
\quad \forall t \in (t_1, t_1 + \ep) ,
   \]
whence~$e^{-t P_1} f_0 = 0$ in~$\Om \setminus \overline{D}_1$ for all~$t \in (0,\ep)$ and finally~$f_0 = 0$ in~$\Om \setminus \overline{D}_1$. But, once more, this is a contradiction and the proof is achieved.
\Fin

\section{Some extensions, comments and open problems}\label{Exten}

   Let us first mention several extensions and variants of 
Theorem~\ref{th.main} and~Corollary~\ref{th.main2}:

\

$\bullet$ First, note that the coefficients~$a_{ij}, b_j, c$ may depend on~$t$.
   The arguments in the proofs hold good if, for instance, they are analytical 
in time.

\

$\bullet$ It is also clear that other boundary conditions on~$\partial\Om\!
\times\!(0,T)$ or~$\partial D\!\times\!(0,T)$ can be imposed in~\eqref{eq.pb}.
   For instance, we can consider the following systems, where an interior input and no boundary input are given:
   \begin{equation}\label{eq.pb-mod}
\begin{cases}
\partial_t u_k + \mathcal{A} u_k = f(x,t) &  \text{in} \  (\Om\!\setminus\!
\overline{D}_k)\times(0,T), 
\\[2mm]
u_k = 0 & \text{on } \  \partial\Om\!\times\!(0,T),
\\[2mm]
\dfrac{\partial u_k}{\partial \nu_A} + m(x) u_k = 0 & \text{on } \  \partial 
D_k\!\times\!(0,T) 
\end{cases}
   \end{equation}
for some nonnegative~$m \in C^0(\Om)$.
   Again, uniqueness can be established arguing as in the proof 
of~Theorem~\ref{th.main}.

\

$\bullet$ We can also prove uniqueness for other right-hand sides~$f$.
   Thus, with~$\mu = \delta_\tau$ (the Dirac distribution at~$\tau$) 
and~$f = f_0(x) \delta_\tau$, it suffices to start from~$u_1$ and~$u_2$ and 
then work with the corresponding primitive functions
   \[
z_k=\int_0^t u_k\, ds , \ \text{ for } \ k=1,2.
   \]

\

$\bullet$ It is possible to establish variants of Theorem~\ref{th.main} 
and~Corollary~\ref{th.main2} with boundary Neumann sources and interior 
observations.
   Thus, for instance, we can start from the systems
   \begin{equation*}
\begin{cases}
\partial_t u_k + \mathcal{A} u_k = 0 & \text{in} \  
(\Om\!\setminus\!\overline{D}_k)\times(0,T), \\[2mm]
\dfrac{\partial u_k}{\partial \nu_A}  = \mu(t)f(x)  & \text{on } 
\  \partial\Om\times(0,T), \\[2mm]
u_k=0 & \text{on } \  \partial D_k\times(0,T)
\end{cases}
   \end{equation*}
  for~$k=1,2$ (where the notation is self-explained).
   By adapting the proof we see that, if for some nonempty open set~$\om\subset
\Om\setminus (\overline{D_1\cup D_2})$ one has~$u_1=u_2$ in~$\om\times(0,T)$, 
then~$D_1 = D_2$.

\

$\bullet$ On the other hand, results similar to Theorem~\ref{th.main} 
and~Corollary~\ref{th.main2} can be proved for more complex linear parabolic 
systems.
   Thus, let us assume that~$(u_k,p_k)$ solves the incomplete quasi-Stokes 
system
   \[
\begin{cases}
\partial_t u_k \!-\! \nu_0\Delta u_k \!+\! (a\cdot\nabla)u_k \!
+\! (u_k\cdot\nabla)b \!+\! \nabla p_k = f(x,t) &  \text{in } 
(\Om\!\setminus\!\overline{D}_k)\!\times\!(0,T) ,
\\[1mm]
\nabla \cdot u_k = 0 &  \text{in} \  (\Om\!\setminus\!\overline{D}_k)\!\times
\!(0,T), 
\\[1mm]
u_k = 0 & \text{on} \ \partial(\Om\!\setminus\!\overline{D}_k)\!\times\!(0,T)
\end{cases}
   \]
for~$k = 1,2$, where the~$D_k$ are as before, $\nu_0 > 0$, $a, b \in 
L^\infty(\Om)^d$ and the components of~$f$ satisfy~\eqref{7a}--\eqref{7c}.
   Let us introduce the notation
   \[
\sigma(u,p) := -p \,\Id + 2\nu_0 \,e(u), \ \text{ where } \ e(u) 
:= \dfrac{1}{2}(\nabla u + (\nabla u)^t)
   \]
and let us assume that
   \begin{equation}\label{eq.coincide}
\sigma(u_1,p_1) \cdot \nu = \sigma(u_2,p_2) \cdot \nu \ \text{ on } 
\ \gamma\!\times\!(0,T),
   \end{equation}
where~$\gamma$ is as above.
   Then, once more, $D_1 = D_2$.
   
   Obviously, we must view here~$u_k$ (resp.~$p_k$) as the velocity field 
(resp.~the pressure distribution) of an incompressible fluid whose particles 
fill the domain~$\Om \setminus \overline{D}_k$ during the time 
interval~$(0,T)$.
   The assumption~\eqref{eq.coincide} means that the {\it normal stresses} 
corresponding to the obstacles~$D_1$ and~$D_2$ coincide on~$\gamma$ at any time.
\

$\bullet$ A result of the same kind can be proved for linearized Boussinesq 
systems.
   More precisely, assume that the triplet~$(u_k,p_k,\theta_k)$ satisfies
   \[
\begin{cases}
\partial_t u_k \!-\! \nu_0\Delta u_k \!+\! (a\!\cdot\!\nabla)u_k \!+\! (u_k\!
\cdot\!\nabla)b \!+\! \nabla p_k \!=\! \theta_k g \!+\! f(x,t) 
& \!\!\!\text{in} \ (\Om\!\setminus\!\overline{D}_k) \!\times\!(0,T) ,
\\[1mm]
\nabla \cdot u_k = 0 & \!\!\!\text{in} \ (\Om\!\setminus\!\overline{D}_k) \!
\times\! (0,T), 
\\[1mm]
\partial_t \theta_k - \kappa_0\Delta \theta_k + a\cdot\nabla\theta_k + u_k\cdot
\nabla q = 0 &  \!\!\!\text{in} \ (\Om\!\setminus\!\overline{D}_k) \!\times\! 
(0,T) ,
\\[1mm]
u_k = 0, \ \ \theta_k = 0 & \!\!\!\text{on} \, \partial(\Om\!\setminus\!
\overline{D}_k) \!\times\! (0,\!T)
\end{cases}
   \]
for~$k = 1,2$ (where~$g \in \mathds{R}^d$, $\kappa_0 > 0$, 
$q \in L^\infty(\Om)$ and the other coefficients are as before) and one has
   \begin{equation}\label{13p}
\sigma(u_1,p_1) \cdot \nu = \sigma(u_2,p_2) \cdot \nu \ \text{ and } \ 
\dfrac{\partial \theta_1}{\partial \nu} = \dfrac{\partial \theta_2}{\partial 
\nu}\ \text{ on } \ \gamma\!\times\!(0,T).
   \end{equation}
   Then~$D_1 = D_2$.
   
   Now, we can interpret that~$u_k$ (resp.~$p_k$ and~$\theta_k$) as the 
velocity field (resp.~the pressure and the temperature distribution) of a 
fluid for which heat effects are important.
   In~\eqref{13p}, we impose that the normal stresses 
and the {\it normal heat fluxes} associated to~$D_1$ and~$D_2$ 
on~$\gamma\!\times\!(0,T)$ are the same. For some related geometric 
inverse problems, see~\cite{ADK-Bous_1}.

\

   Also, some comments on~Theorem~\ref{th.main3} may be given:

\

$\bullet$ Again, as in~\eqref{eq.pb-mod}, other boundary conditions 
on~$\partial\Om\!\times\!(0,T)$ or~$\partial D\!\times\!(0,T)$ can be imposed.

\

$\bullet$ Results similar to Theorem~\ref{th.main3} can also be established 
for the incomplete quasi-Stokes and Boussinesq systems corresponding to 
different viscosities, heat diffusions and coefficients.

   For example, the partial differential equations can be of the form
   \[
\begin{cases}
\partial_t u_k - \nu_k\Delta u_k + (a_k\cdot\nabla)u_k + y_k\!\times\!u_k 
+ \nabla p_k = f(x,t) ,
\\[1mm]
\nabla \cdot u_k = 0 ,
\end{cases}
   \]
where the~$\nu_k > 0$, the~$a_k, y_k \in \R^d$ and the components of~$f$ 
satisfy~\eqref{7a}, \eqref{7b} and~\eqref{7d}.
   If one has
   \[
u_1 = u_2 \ \text{ on } \om \times (0,T), 
   \]
then necessarily~$D_1 = D_2$.

   Here, an interesting interpretation can be given:
   assume that the oceans of two different planets of the same size that 
move with different translational and rotational 
velocities are modelled by quasi-Stokes systems with different viscosities;
   also, assume that the corresponding ocean velocity fields are identical 
in a small region at any time;
   then, the shapes of the associated continents necessarily coincide.

Together with these extensions, several open questions are in order:
   
   \
   
$\bullet$ What about multiply connected~$D_1$ and~$D_2$?
   
   \
   
$\bullet$  Can we also get uniqueness for more general inputs of the form
   \[
F = F(x,t) \ \text{ with } \ \text{Supp\,} F(\cdot, t)\subset \Om\!\setminus
\!\overline{D_1 \cup D_2} \ \text{ for all } \ t?
   \]

 \

$\bullet$ Note that we can get uniqueness for the quasi-Stokes systems 
considered above by assuming that just~$d-1$ components of~$\sigma(u_1,p_1) 
\cdot \nu$ and~$\sigma(u_2,p_2) \cdot \nu$ coincide on~$\gamma\times(0,T)$.
   The reason is that this, together with the incompressibility of~$u_1$ 
and~$u_2$ and the requirement~$u_1 = u_2$ on~$\partial\Om\times(0,T)$, 
suffices to yield~$u_1 = u_2$ in~$(\Om\!\setminus\!\overline{D_1 \cup D_2})
\times(0,T)$, see~\cite{JLL-EZ};
   see also~\cite{Coron-SGR, ADK-Bous_1, EFC-N-1}.

However, for~$d= 3$, it is unknown whether the result holds when only one 
component coincides (in this case, there are domains for which unique 
continuation fails, although it does not seem easy to characterize 
them completely, see~\cite{JLL-EZ}).

Similar considerations apply for linearized Boussinesq systems:
if we assume that~$d-1$ components of~$\sigma(u_1,p_1) \cdot \nu$ 
and~$\sigma(u_2,p_2) \cdot \nu$ and the normal derivatives of~$\theta_1$ 
and~$\theta_2$ coincide on~$\gamma\!\times\!(0,T)$, uniqueness holds;
   nevertheless, if~$d = 3$ and only one component of the normal stresses and 
the normal heat fluxes coincide, uniqueness is open.

\

$\bullet$   As a final additional system, let us consider the linearized at 
zero Oldroyd model for a visco-elastic fluid and a related open question; 
see~\cite{ADK-EFC, EFC-FGG-RRO, Joseph} for a motivation. Thus, let 
the~$(u_k,p_k,\tau_k)$ satisfy
   \[
\begin{cases}
\partial_t u_k - \nu \Delta u_k + \nabla p_k = \nabla \cdot \tau_k + f(x,t) 
&  \text{in} \  (\Om\!\setminus\!\overline{D}_k)\times(0,T) ,
\\[1mm]
\nabla \cdot u_k = 0 &  \text{in} \  (\Om\!\setminus\!\overline{D}_k)
\times(0,T) ,
\\[1mm]
\partial_t \tau_k + \alpha\,\tau_k =  2\beta\,e(u_k) &  \text{in} \  
(\Om\!\setminus\!\overline{D}_k)\times(0,T) ,
\\[1mm]
u_k = 0 & \text{on } \  \partial(\Om\!\setminus\!\overline{D}_k) \times(0,T)
\end{cases}
   \]
for~$k=1,2$, where~$\alpha, \beta > 0$ and~$u_k$, $p_k$ and~$e(u_k)$ are 
interpreted as in the above quasi-Stokes system. Here, the~$\tau_k$ are the 
so called {\it extra-elastic tensors.}
   They must be viewed as tensor stresses originated by elastic (internal) 
forces exerted by particles.

   Assume that we have~\eqref{eq.coincide}.
   Then, do we have~$D_1 = D_2$?

\subsection*{Acknowledgements}
   The first author was partially supported by Grant~PID2021-126813NB-I00, funded by MCIN/AEI/10.13039/501100011033 and by ``ERDF A way of making Europe'' and Grant~IT1615-22, funded by the Basque Government (Spain).
   On the other hand, the second and third authors were partially supported by Grants~PID2020-114976GB-I00 and~PID2024-158206NB-I00, funded by MICIU/AEI (Spain).
   The fourth author was supported by Grant-in-Aid for Scientific Research (A) 20H00117 and Grant-in-Aid for Challenging Research (Pioneering) 21K18142 of Japan Society for the Promotion of Science.

\end{document}